\documentclass[12pt]{amsart}

\usepackage{amsmath, amsthm, amscd, amsfonts, amssymb, enumerate,
  verbatim, newlfont, calc, graphicx, color}
\usepackage[bookmarksnumbered, colorlinks, plainpages]{hyperref}
\usepackage{tikz} 
\usepackage{doi}
\usepackage{cancel}
 
\newtheorem{theorem}{Theorem}[section]       
\newtheorem{question}[theorem]{Question}    
 
\newtheorem{lemma}[theorem]{Lemma}
\newtheorem{proposition}[theorem]{Proposition}

\newtheorem{corollary}[theorem]{Corollary}
\newtheorem{definition}[theorem]{Definition}

\newtheorem{exercise}[theorem]{Exercise}

\usepackage{mathtools}

\usepackage{pstricks}
\usepackage{epsfig}
\usepackage{pst-grad}
\usepackage{pst-plot}
\usepackage{subfig}

\begin{document}

\author[ W. Hochst\"attler   and  M. Nasernejad]{Winfried  Hochst\"attler$^{1}$  and   Mehrdad ~Nasernejad$^{2,*}$}

\title[A classification of Mengerian $4$-uniform hypergraphs] {A classification of Mengerian $4$-uniform hypergraphs derived from graphs}
 


\thanks{$^*$Corresponding author}

\thanks{E-mail addresses:    winfried.hochstaettler@fernuni-hagen.de  and  m$\_$nasernejad@yahoo.com} 
\maketitle

\begin{center}
{\it

$^{1}$FernUniversit\"{a}t in Hagen, Fakult\"{a}t f\"{u}r Mathematik und Informatik,  58084 Hagen,  Germany\\
$^{2}$Univ. Artois, UR 2462, Laboratoire de Math\'{e}matique de Lens \\ (LML),   F-62300 Lens, France

}
\end{center}

\vspace{0.4cm}

%%% ----------------------------------------------------------------------

\begin{abstract}
In this paper, we give a classification of all Mengerian $4$-uniform hypergraphs derived from graphs.  
\end{abstract}

%%% ----------------------------------------------------------------------

\vspace{0.4cm}

%%%%%%%%%%%%%%%%%%%%%%%%%%%%%%%%%%%%%%%%%%%%%%%%%%%%%%%%%%%%
%%%%%%%%%%%%%%%%%%%%%%%%%%%%%%%%%%%%%%%%%%%%%%%%%%%%%%%%%%%%
%%%%%%%%%%%%%%%%%%%%%%%%%%%%%%%%%%%%%%%%%%%%%%%%%%%%%%%%%%%%

\section{Introduction}

Assume that   $R$ is   a commutative Noetherian ring and $I$ is   an ideal of $R$. Then a  prime ideal $\mathfrak{p}\subset  R$ is an {\it associated prime} 
    of $I$ if there exists an element $h$ in $R$ such that $\mathfrak{p}=(I:_R h)$, where $(I:_Rh)=\{r\in R \mid rh\in I\}$.
     The  {\it set of associated primes} of $I$, denoted by  $\mathrm{Ass}_R(I)$, is the set of all prime ideals associated to  $I$.  
      Brodmann \cite{BR} showed that the sequence $\{\mathrm{Ass}_R(I^k)\}_{k \geq 1}$ of associated prime ideals is stationary  for large $k$. 
     This means that  there exists a positive integer $k_0$ such that $\mathrm{Ass}_R(I^k)=\mathrm{Ass}_R(I^{k_0})$ for all $k\geq k_0$.  The  minimal such $k_0$ is called the {\it index of stability}   of  $I$ and $\mathrm{Ass}_R(I^{k_0})$ is called the {\it stable set } 
 of associated prime ideals of  $I$, which is denoted by $\mathrm{Ass}^{\infty }(I).$ 
 
Several  questions arise  in the   context of Brodmann’s result. 
 An ideal $I\subset R$ is said to be  {\it normally torsion-free} if, for all $k\geq 1$, we have  $\mathrm{Ass}_R(I^k)\subseteq \mathrm{Ass}_R(I)$, 
  see \cite[Definition 4.3.28]{V1}. In particular, if $I$ is a square-free monomial ideal in a polynomial ring $R=K[x_1, \ldots, x_n]$ over a field $K$, then 
  $I$ is  normally torsion-free  if, for all $k\geq 1$, we have  $\mathrm{Ass}_R(I^k)=\mathrm{Ass}_R(I)$.
    
   Generally, finding classes of normally torsion-free ideals have been a theme  of many  papers, nevertheless, some 
  classes of such  ideals emanate from graph theory. Particularly, it is well-known that a finite simple undirected graph is bipartite if and only if its edge ideal is normally torsion-free, if and only if its cover ideal is normally torsion-free, see \cite{GRV, SVV}  for more information. 
Furthermore,  in \cite[Theorem  3.3]{N3}, the authors stated   that  the Alexander dual of the monomial
ideal generated by the paths of maximal lengths in  a rooted starlike tree is normally torsion-free. 
In  \cite[Theorem 3.2]{KHN1}, it has been shown  that   the Alexander dual of path ideals generated by all paths of length $2$ in rooted trees are normally
torsion-free. 

Let $G$ be an undirected and simple graph. More recently, it has been verified in \cite[Theorem 2.8]{NQ} that if $G$ is  a strongly chordal graph, then  both $NI(G)$ and $DI(G)$ are  normally torsion-free, where $NI(G)$ (respectively, $DI(G)$) stands for the closed neighborhood ideal of $G$ (respectively, dominating ideal of $G$). 
Moreover, when $C_n$ denotes the simple cycle graph on $n$ vertices, then based on  \cite[Theorem 4.3]{NQ},    $DI(C_n)$ is normally torsion-free if and only if $n \in \{3,6,9\}$, 
the definitions of closed neighborhood ideals and  dominating ideals and their properties can be found in \cite{NBR, NQBM}. 
 
In addition,   Lemma 5.15  in \cite{NQKR} gives a  characterization of  all normally torsion-free $t$-spread principal Borel ideals. In fact, 
 a monomial $x_{i_1} x_{i_2} \cdots x_{i_d} \in R=K[x_1, \ldots, x_n]$ with $i_1 \leq i_2 \leq \cdots  \leq i_d$ is called {\it $t$-spread} if $i_j -i_{j-1} \geq t$ for all $j=2, \ldots, d$.  A monomial ideal $I$ in $R$ is called a {\it $t$-spread monomial ideal} if it is generated by $t$-spread monomials. Also,  $I$ is called a {\it $t$-spread strongly stable ideal} if for all $t$-spread monomials $u\in \mathcal{G}(I)$, all $j\in \mathrm{supp}(u)$
 and all $1\leq i <j$ such that $x_i(u/x_j)$  is a $t$-spread monomial, it follows that $x_i(u/x_j)\in I$. 
Moreover, $I$  is called a {\it $t$-spread principal Borel} if there exists a $t$-spread monomial $u\in \mathcal{G}(I)$ such that  $I$  is the smallest $t$-spread strongly stable ideal which contains $u$; we denote it as $I=B_t(u)$. By \cite[Lemma 5.15(i)]{NQKR}, if  $u=x_ax_bx_n$ is a $t$-spread monomial, 
$I=B_t(u)$, and $b < 2t+1$, then $I$  is normally torsion-free.
 In this direction,  in  \cite[Theorem 4.6]{MNQ}, the authors showed  that the edge ideal of any  $\mathbf{t}$-spread $d$-partite hypergraph   
 $K^{\textbf{t}}_{V}$ is normally torsion-free. More  information about $t$-spread monomial ideals and  vector-spread monomial ideals can be found 
 in \cite{EHQ, F}. 
 
 Furthermore, Herzog et. al., in \cite[Corollary 4.6]{HRV},  proved  that every  transversal polymatroidal ideal  is normally torsion-free.  In particular, Olteanu \cite{OL}  characterized all the lexsegment ideals which are  normally torsion-free non-square-free monomial ideals.

 In the following text, we focus on normally torsion-free square-free monomial ideals. Indeed, one  of the motivations of this work originates from this fact that
    normally torsion-free square-free monomial ideals are closely related to Mengerian hypergraphs.  A hypergraph  is called 
\textit{Mengerian} if it satisfies a certain min-max equation, which is known as the Mengerian property in hypergraph theory or has the max-flow min-cut property in integer programming. In fact,  it is well-known that if $\mathcal{C}$ is  a clutter, i.e. an antichain w.r.t. the partial order
induced by inclusion, and $I=I(\mathcal{C})$  its edge ideal, then $I$ is  normally torsion-free if and only if  $I^{k}=I^{(k)}$ for all $k\geq 1$, where $I^{(k)}=\bigcap_{\mathfrak{p}\in \mathrm{Min}(I)}(I^kR_\mathfrak{p}\cap R)$ denotes the  $k$-th symbolic  power  of $I$,  if and only if $\mathcal{C}$ is Mengerian, if and only if
 $\mathcal{C}$  has the max-flow min-cut  (MFMC for short)  property,  see  \cite[Theorem 14.3.6]{V1} for more details. 
 
So far,  there exist two well-known  classifications of Mengerian $2$-uniform hypergraphs and $3$-uniform hypergraphs derived from graphs as follows:

\begin{itemize}
\item[(i)] (\cite[Theorem  5.9]{SVV})  Let  $G$  be a simple undirected  graph and $I$ its edge ideal. Then $G$ is bipartite if and only if $I$ is normally torsion-free. 
\item[(ii)] (\cite[Theorem 3.9]{AB}) Let $G$ be a connected  simple undirected  graph and  $t \geq 2$  be an integer, and let   $J=I_3(G)$ 
be the cubic path ideal of  $G$. Then, $J^{(n)} = J^n$ for all  $n \geq  1$  if and only if $G$  is a path graph $P_t$ 
 or $G$ is the cycle $C_{3k}$ when $k = 1, 2, 3$.
\end{itemize}

In addition, according to \cite[Proposition 3.2]{NQ}, the $t$-path ideals of path graphs are normally torsion-free  for all $t\geq 0$. 

The main aim of this work is to provide a classification of  all Mengerian $4$-uniform hypergraphs derived from graphs as we will give  in the following theorem:

\textbf{Theorem  \ref{Main-Result}.} 
If $G=(V, E)$ is a connected graph, then $\mathcal{H}_3(G)$ is Mengerian if and only if $|V|=4$, or $G=C_8$, or $G$ is a path with double stars, or 
a star with an additional edge  between two of its leaves.

\bigskip
It is possible to ask whether the techniques used in this paper can be generalized to the general case? Our computations and observations  show that investigating the cases greater than $4$ is too complicated, and we may need to change our technique as in \cite{SVV} the authors used the Rees Algebra and classic methods in commutative algebra  for proving the case $2$, and in \cite{AB} the authors utilized the symbolic powers technique for showing the case $3$, and the present authors used the TDI (totally dual integral)    method for proving the case $4$;  however, we left an open question in this direction at the end of this paper, see Question \ref{Open Question}.

\bigskip 
Throughout this paper, all graphs  are simple, finite, connected, and undirected.  Also, any  necessary definitions related to graph theory  (respectively, monomial ideals) can be found in  \cite{West} (respectively, \cite{HH1}).

%%%%%%%%%%%%%%%%%%%%%%%%%%%%%%%%%%%%%%%%%%%%%%%%%%%%%%%%%%%%
%%%%%%%%%%%%%%%%%%%%%%%%%%%%%%%%%%%%%%%%%%%%%%%%%%%%%%%%%%%%
%%%%%%%%%%%%%%%%%%%%%%%%%%%%%%%%%%%%%%%%%%%%%%%%%%%%%%%%%%%%

\section{Preliminaries}

In what follows, we  recollect  some necessary definitions related to  hypergraph theory, which can be found in \cite{Berge, V1}. 

A finite {\it hypergraph} $\mathcal{H}$ on a vertex set $V({\mathcal{H}})=\{x_1,x_2,\ldots,x_n\}$ is a collection of  hyperedges  $E({\mathcal{H}})=\{ E_1, \ldots, E_m\}$ with $E_i \subseteq V({\mathcal{H}})$
 for all $i=1, \ldots,m$.  The {\em incidence matrix} of $\mathcal{H}$ is an $m \times n$ matrix $A=(a_{ij})$ such that $a_{ij}=1$ if $x_j \in E_i$, and $a_{ij}=0$ otherwise. 
A hypergraph $\mathcal{H}$ is called {\it simple}  if $E_i \subseteq  E_j$ implies  $i = j$. Simple hypergraphs are also known as {\em clutters}. Moreover, if $|E_i|=d$  for all $i=1, \ldots, m$, then $\mathcal{H}$ is called a {\em $d$-uniform} hypergraph. A $2$-uniform hypergraph $\mathcal{H}$ is just a finite simple graph.

Let $V(\mathcal{H})=\{x_1,x_2,\ldots, x_n\}$ and $R=K[x_1, \ldots, x_n]$ be the polynomial ring in $n$ variables  over a field $K$. The {\it edge ideal} of $\mathcal{H}$ is given by
$$I(\mathcal{H}) = \left(\prod_{x_j\in E_i} x_j : E_i\in  E({\mathcal{H}})\right).$$

Let $G$ be a simple finite connected  undirected graph. The {\em $t$-path hypergraph} of $G$, denoted by $\mathcal{H}_t(G)$, is a $(t+1)$-uniform hypergraph on the vertex set $V(G)$ such that $e=\{v_1, \ldots, v_{t+1}\}$ is an edge if contains  a  path of length $t$ in $G$.  The edge ideal of $\mathcal{H}_t(G)$ is called the {\em $t$-path ideal} of $G$ and is denoted by $I_t(G)$.  When $t=1$, then $I_1(G)$ is simply the {\em edge ideal}  $I(G)$ of $G$. If there is no $t$-path in $G$, then we set $I_t(G)=0$.  

\bigskip
To understand Theorem \ref{Villarreal1}, we require to review the following definitions from combinatorial optimization. 
\begin{definition} (\cite[Definition 14.3.4]{V1})
\em{
 A clutter $\mathcal{C}$ satisfies the \textit{max-flow min-cut (MFMC)} property  if both sides of the LP-duality equation
\begin{equation} \label{MFMC}
 \min\{\langle{c}, x\rangle | x\geq 0; Ax\geq \textbf{1}\}=\max\{\langle y, \textbf{1}\rangle | y\geq 0; y^\top A\leq c^\top\}
\end{equation}
 have integral optimum solutions $x$ and $y$ for each nonnegative integral vector $c$. 
  The system $x \geq  0; Ax \geq  1$ is called  \textit{totally dual integral (TDI)} if the
maximum in Eq. (\ref{MFMC}) has an integral optimum solution $y$ for each integral
vector $c$ with finite maximum.
}
\end{definition}

For a hypergraph $\mathcal{H}$, the   {\it deletion} $\mathcal{H}\setminus x$ at a vertex $x\in V(\mathcal{H})$ has the effect of setting $x=0$ in the edge ideal $I:=I(\mathcal{H})$, and the {\it contraction} $\mathcal{H}/x$ has the effect of setting $x=1$ in the edge ideal $I(\mathcal{H})$, see 
\cite[Definition 6.5.2]{V1}. 
 We call an ideal $I'$ a {\it minor}   of  a square-free monomial ideal $I$ if $I'$ can be obtained from $I$ by a sequence of taking quotients (=deletion) and localizations (=contraction) at the variables, refer to   \cite[Definition  6.5.3]{V1}.

 A vertex $x\in V(\mathcal{H})$  is called an {\it isolated vertex}  if $\{x\}\in E(\mathcal{H})$. It follows 
from the definition that if $x$ is an isolated vertex of $\mathcal{H}$, then $\{x\}$  is the only edge in  $\mathcal{H}$ 
 that contains $x$. A {\it vertex cover}  of  $\mathcal{H}$ is a set of vertices that has a  non-empty
intersection with all of the edges of $\mathcal{H}$.  The minimum cardinality of a vertex cover of $\mathcal{H}$
is denoted by $\tau(\mathcal{H})$, consult \cite[Page 53]{Berge}. It should be noted that  by the correspondence between minimal primes 
 and minimal vertex covers, $\tau(I)$ is also the height of $I$  (cf. \cite[Proposition 13.1.6]{V1}).

We say the  generators of a square-free monomial ideal $I$ as being {\it independent} or a {\it matching} if the corresponding edges of the
associated hypergraph are pairwise disjoint; that is, the generators have disjoint supports. We  denote the maximum
cardinality of an independent set in $I$ by $\nu(I)$, see \cite[Page 64]{Berge}. Furthermore, note that a subset  of the monomial generators of  $I$
is independent if and only if it forms a regular sequence. Thus,  $\nu(I)$ is equal to the \textit{monomial grade}  of $I$, denoted by 
$\mathfrak{m}$-$\mathrm{grade}(I)$,  where the monomial grade of an ideal $I$ is the maximum length of a regular sequence of monomials in $I$, 
see \cite[Definition 13.1.5 and Proposition 13.1.6]{V1}. By weak linear programming duality we always have  $\tau(I) \geq \nu(I)$.

A square-free monomial ideal $I$  is said to satisfy the {\it K\"onig property }
if  $\tau(I) = \nu(I)$, refer to  \cite[Definition 13.1.4]{V1}.   Consequently,  $I$ satisfies the K\"onig property if and only if  
$\mathrm{grade}(I)=\mathrm{ht}(I)=\mathfrak{m}$-$\mathrm{grade}(I)$. Moreover, a  square-free monomial ideal $I$ has the {\it packing property}     if  $I$  and all of its proper minors satisfy the K\"onig property (cf. \cite[Definition 14.3.16]{V1}). 

\begin{definition} (\cite[Definition 14.3.5]{V1}) 
\em{
 A clutter  $\mathcal{C}$ is called \textit{Mengerian} if $\nu({\mathcal{C}}^a)
  = \tau({\mathcal{C}}^a)$ for all $a \in \mathbb{N}^n$, where
  $\mathbb{N}$ denotes the set of nonnegative integers. In other   words, $\mathcal{C}$ is Mengerian if all its multiplications of
  vertices have the K\"onig property.
  }
\end{definition}

 In 1990 \cite{CC}, Michele  Conforti and  G\'e{r}ard   Cornu\'e{j}ols     made up  the following conjecture, and it is still open.  

\textbf{Conjecture.} (Conforti-Cornu\'e{j}ols  conjecture) 
  A hypergraph has the packing property if and only if it is Mengerian.

We finally  summarize all of  the above notions in the following theorem. 

\begin{theorem} (\cite[Theorem 14.3.6]{V1})  \label{Villarreal1}
 Let $\mathcal{C}$  be a clutter and let  $A$ be its incidence matrix. The following are equivalent:
 \begin{itemize}
 \item[(i)] $\mathrm{gr}_I(R)$ is reduced, where $I=I(\mathcal{C})$ is the edge ideal of $\mathcal{C}$. 
 \item[(ii)] $R[It]$ is normal  and $\mathcal{Q}(A)$ is an integral polyhedron. 
 \item[(iii)]  $x\geq 0$; $Ax\geq \textbf{1}$   is a $\mathrm{TDI}$ system. 
 \item[(iv)]  $\mathcal{C}$  has the max-flow min-cut (MFMC) property.
 \item[(v)]  $I^i = I^{(i)}$ for $i \geq 1$. 
 \item[(vi)]  $I$ is normally torsion-free, i.e.,  $\mathrm{Ass}(I^i) \subseteq  \mathrm{Ass}(I)$ for $i\geq 1$. 
 \item[(vii)]  $\mathcal{C}$  is Mengerian, i.e., $\nu({\mathcal{C}}^a) = \tau({\mathcal{C}}^a)$ for all  $a \in \mathbb{N}^n$.
 \end{itemize}
\end{theorem}

%++++++++++++++++++++++++++++++++++++++++++++++++++++++++++++++++++++++++++++++++++
%++++++++++++++++++++++++++++++++++++++++++++++++++++++++++++++++++++++++++++++++++

%%%%%%%%%%%%%%%%%%%%%%%%%%%%%%%%%%%%%%%%%%%%%%%%%%%%%%%%%%%%
%%%%%%%%%%%%%%%%%%%%%%%%%%%%%%%%%%%%%%%%%%%%%%%%%%%%%%%%%%%%
%%%%%%%%%%%%%%%%%%%%%%%%%%%%%%%%%%%%%%%%%%%%%%%%%%%%%%%%%%%%

\section{Main result}

 Given a connected graph $G=(V, E)$ with $|V|\geq 4$, we consider the $4$-uniform 
 hypergraph $\mathcal{H}_3(G)=(\mathcal{V}, \mathcal{F})$ such that
 $e=\{v_1, v_2, v_3, v_4\} \in \mathcal{F}$ if and only if contains a
 path of length $3$ in $G$.  Let $A$ denote the incidence matrix of
 $\mathcal{H}_3(G)$.  Then $\mathcal{H}_3(G)$ is Mengerian if and only
 if the system $Ax\geq \mathbf{1}$ with $x\geq 0$ is TDI. This in
 particular implies that the system is integral. In this case, we call
 the hypergraph \textit{ideal}. It should be noted that if a
 hypergraph is not ideal, then it is not Mengerian.
Furthermore, every empty clutter is considered Mengerian.

To show Proposition \ref{Pro.1} here below, we require to employ the following definition, exercise, and theorem. Given a  monomial ideal $I$, we denote by  $\mu(I)$ the minimal number 
of  monomial generators of $I$. 

%%%%%%%%%%%%%%%%%%%%%%%%%%%%%%%%%%%%%%%%%%%%%%%%%%%%%%%%%

\begin{definition} (\cite[Definition 4.3.22]{V1}) \label{Def.4.3.22}
\em 
{
 Let $I$ be an ideal of a ring $R$ and $\mathfrak{p}_1, \ldots, \mathfrak{p}_r$  the minimal primes of $I$. Given an integer 
 $n \geq 1$,  the $n$th \textit{symbolic power} of $I$ is defined to be the ideal
 $I^{(n)} = \mathfrak{q}_1 \cap \cdots \cap \mathfrak{q}_r$,  where $\mathfrak{q}_i$  is the primary component of $I^n$ corresponding to $\mathfrak{p}_i$. 
 }
 \end{definition}
 
 %%%%%%%%%%%%%%%%%%%%%%%%%%%%%%%%%%%%%%%%%%%%%%%%%%%%%%%%%%%%
 
 \begin{exercise} (\cite[Exercise 6.1.25]{V1})\label{Ex.6.1.25}
 If  $\mathfrak{q}_1, \ldots, \mathfrak{q}_r$  are primary monomial ideals of $R$ with non-comparable radicals and $I$ 
is an ideal such that $I =\mathfrak{q}_1 \cap  \cdots \cap  \mathfrak{q}_r$, then $I^{(n)} =\mathfrak{q}^n_1 \cap  \cdots \cap  \mathfrak{q}^n_r$. 
\end{exercise}

%%%%%%%%%%%%%%%%%%%%%%%%%%%%%%%%%%%%%%%%%%%%%%%%%%%%%%%%%%%%

\begin{theorem} (\cite[Theorem 4.8]{MN}) \label{Montana-Betancourt}
Let $I$ be a square-free monomial ideal. If $I^n=I^{(n)}$ for every $n\leq \lceil\frac{\mu(I)}{2}\rceil$, then $I^n=I^{(n)}$ for every 
$n\in \mathbb{Z}_{>0}$. 
\end{theorem}

%%%%%%%%%%%%%%%%%%%%%%%%%%%%%%%%%%%%%%%%%%%%%%%%%%%%%%%%%%%%

\begin{proposition} \label{Pro.1}
Let $C_8=(V(C_8), E(C_8))$ denote the cycle graph with vertex set $V(C_8)=\{x_1, \ldots, x_8\}$ and the following 
 edge set $$E(C_8)=\{\{x_i, x_{i+1}\} : i=1, \ldots, 7\} \cup \{\{x_1, x_8\}\}.$$ 
 Let $\mathcal{H}_3(C_8)$ be the $4$-uniform hypergraph on the vertex set  $V(C_8)$.  Then $\mathcal{H}_3(C_8)$ is Mengerian. 
\end{proposition}

\begin{proof} 
 Let  $I_3(C_8)$ be  the edge ideal of $\mathcal{H}_3(C_8)$, and $R=K[x_1, \ldots, x_8]$.
For convenience of notation, we put $J:=I_3(C_8)$. It is routine   to  check  that 
\begin{align*}
 J=&(x_1x_2x_3x_4, x_2x_3x_4x_5, x_3x_4x_5x_6, x_4x_5x_6x_7, x_5x_6x_7x_8, \\
 &x_6x_7x_8x_1, x_7x_8x_1x_2, x_8x_1x_2x_3).
 \end{align*}
Hence, we deduce that $\mu(J)=8$. Using {\it Macaulay2} \cite{GS}, we obtain
\begin{align*}
\mathrm{Ass}(J)=&\mathrm{Ass}(J^2)=\mathrm{Ass}(J^3)=\mathrm{Ass}(J^4)\\
=&\{(x_5,x_1), (x_6,x_2),  (x_7,x_3),  (x_8,x_4),   (x_6,x_3,x_1),  (x_6,x_4,x_1),\\
 & (x_7,x_4,x_1), (x_7,x_4,x_2), (x_7,x_5,x_2),  (x_8,x_5,x_2),  (x_8,x_5,x_3), \\
 & (x_8,x_6,x_3)\}. 
\end{align*}
Since $J$ is a square-free monomial ideal, we  deduce from \cite[Corollary 1.3.6]{HH1} that $\mathrm{Ass}(J)=\mathrm{Min}(J)$, in particular, 
$J=\bigcap_{\mathfrak{p}\in \mathrm{Min}(J)}\mathfrak{p}$. 
Hence, we have $\mathrm{Ass}(J^n)=\mathrm{Min}(J)$ for all $1\leq n \leq 4$. This means that  $J^n$ (for all $1\leq n \leq 4$) has no embedded associated prime ideals. We can derive  from Definition \ref{Def.4.3.22} and Exercise \ref{Ex.6.1.25}
  that $J^n=J^{(n)}=\bigcap_{\mathfrak{p}\in \mathrm{Min}(J)}\mathfrak{p}^n$ for all $1\leq n \leq 4$.  It follows  now from  Theorem \ref{Montana-Betancourt} that  $J^n=J^{(n)}$ for any $n\geq 1$. One can immediately conclude from Theorem \ref{Villarreal1} that $J$ is normally torsion-free, and so 
$\mathcal{H}_3(C_8)$ is Mengerian, as required. 
\end{proof}

%%%%%%%%%%%%%%%%%%%%%%%%%%%%%%%%%%%%%%%%%%%%%%%%%%%%%%%%%%%%
%%%%%%%%%%%%%%%%%%%%%%%%%%%%%%%%%%%%%%%%%%%%%%%%%%%%%%%%%%%%

 To formulate Lemma \ref{Lem.01}, we need the following definitions and corollary.  
 
 \begin{definition}
 \em{
 A tree $T$ is called a \textit{path with double stars} if it has at most two vertices which are adjacent to leaves. Note that this includes stars and paths. 
 }
 \end{definition}

\begin{definition} (\cite[Definition 1.5.5]{V1})
\em{
 A matrix $A$ is called \textit{totally unimodular} if each $i\times i$ subdeterminant of $A$ is $0$ or $\pm 1$ for all $i \geq 1$. 
 }
 \end{definition}
 
\begin{corollary} (\cite[Corollary 14.3.14]{V1}) \label{Cor.1}
 If $A$ is totally unimodular, then $\mathrm{gr}_I(R)$ is reduced, and so the clutter is    Mengerian.
 \end{corollary} 

%%%%%%%%%%%%%%%%%%%%%%%%%%%%%%%%%%%%%%%%%%%%%%%%%%%%%%%%%%%%
%%%%%%%%%%%%%%%%%%%%%%%%%%%%%%%%%%%%%%%%%%%%%%%%%%%%%%%%%%%%

 \begin{lemma} \label{Lem.01}
  If $G=(V(G), E(G))$ is a connected graph and $A$ is  the   incidence matrix of    $\mathcal{H}_3(G)$ 
  such that  $G$ is a path with double stars, where the path has at least two edges, or $G$ is a star plus an edge    between two of its leaves, then  $A$ is totally
   unimodular,  and thus  $\mathcal{H}_3(G)$ is   Mengerian.
\end{lemma}

\begin{proof}
  We proceed by induction on the number of vertices. If $|V(G)|\le 4$,  then the 
  hypergraph is either empty  or consists of one hyperedge and the
  matrix is one all $1$-row. In both cases, the matrix  $A$ is
  totally unimodular. Let  $|V(G)|>4$ and   $v$ be a vertex of degree $1$ and $A'$ be a square submatrix of 
  $A$. If $A'$ does not intersect the column corresponding to $v$, then 
  $\det(A') \in \{0, \pm 1\}$ by inductive assumption applied to $G
  \setminus v$. If $A'$ intersects the column of $A$, then this column
  can have at most one $1$. If it is all zero, then  we have $\det(A') = 0$. Otherwise, 
     using Laplacian expansion and again inductive assumption applied to $G
  \setminus v$,  we again find   $\det(A') \in \{0, \pm 1\}$.  This completes the inductive step,
   and so the claim has been shown by induction.  The last assertion can be deduced  according to  Corollary \ref{Cor.1}. 
\end{proof}

%%%%%%%%%%%%%%%%%%%%%%%%%%%%%%%%%%%%%%%%%%%%%%%%%%%%%%%%%%%%%
%%%%%%%%%%%%%%%%%%%%%%%%%%%%%%%%%%%%%%%%%%%%%%%%%%%%%%%%%%%%%
As was pointed out by an anonymous referee, the above proof does not
apply to the case where we have a path consisting of just one edge
with double stars. Thus, we need an extra lemma.

 \begin{lemma} \label{Lem.01a}
  If $G=(V(G), E(G))$  is a path with double stars, where the path has exactly one  edge, then  $A$ is totally
   unimodular,  and thus  $\mathcal{H}_3(G)$ is   Mengerian.
\end{lemma}

\begin{proof}
  We construct a network from $G$ as follows. Subdivide the inner
  edge, and direct all edges of the resulting tree such that all paths using four edges are
  directed. Call that directed tree $T$.  Now add all arcs from all
  leaves of the first star to the leaves of the other star, i.e.\ the
  corresponding complete bipartite digraph, resulting in a directed graph $D$.  Consider the
  totally unimodular {\em network matrix} (see~\cite{Schrijver86}
  Chapter 19, Example 4) $\tilde A$ we get from the digraph $D$ and its tree $T$. 

  Then $A$ is the submatrix of $\tilde A$ we get by deleting all
  columns corresponding to arcs in $T$. As a submatrix of a totally
  unimodular matrix, $A$ is totally unimodular as well.
\end{proof}

%%%%%%%%%%%%%%%%%%%%%%%%%%%%%%%%%%%%%%%%%%%%%%%%%%%%%%%%%%%%%
%%%%%%%%%%%%%%%%%%%%%%%%%%%%%%%%%%%%%%%%%%%%%%%%%%%%%%%%%%%%%

\begin{proposition} \label{Pro.2}
If  $G=(V(G), E(G))$ is a connected graph such that $|V(G)|=4$, or $G=C_8$, or $G$ is a path with double stars, or $G$ is  a star plus an edge between two of its leaves, then $\mathcal{H}_3(G)$ is Mengerian. 
\end{proposition}

 \begin{proof}
   The claim clearly is true  if $|V(G)|=4$. If $G=C_8$, by virtue of 
   Proposition \ref{Pro.1}, we can conclude that  $\mathcal{H}_3(C_8)$ is Mengerian.
   In any other case,  $\mathcal{H}_3(G)$ is Mengerian on account of  Lemmata~\ref{Lem.01} and ~\ref{Lem.01a}.
 \end{proof}
 
 %%%%%%%%%%%%%%%%%%%%%%%%%%%%%%%%%%%%%%%%%%%%%%%%%%%%%%%%%%%%
%%%%%%%%%%%%%%%%%%%%%%%%%%%%%%%%%%%%%%%%%%%%%%%%%%%%%%%%%%%%

 In the remaining cases,  we will  show that the hypergraphs are not
 Mengerian because they are not ideal. For this  purpose,  we present in
 each case a fractional vertex of the covering polyhedron. We use the following well-known fact from linear programming.

 \begin{lemma}\label{Lem.02}(see \cite[Theorem 8.4 and (23) in Section 8.5]{Schrijver86})
   Let $P=\{Ax \le b\}$ denote a polyhedron. Then $v \in P$ is a
   vertex of $P$ if and only if there is a subsystem $A'x \le b'$ of
   $Ax\le b$ such that $A'$ has full column rank and $A'v=b'$.
 \end{lemma}

%%%%%%%%%%%%%%%%%%%%%%%%%%%%%%%%%%%%%%%%%%%%%%%%%%%%%%%%%%%%
%%%%%%%%%%%%%%%%%%%%%%%%%%%%%%%%%%%%%%%%%%%%%%%%%%%%%%%%%%%%

 \begin{lemma} \label{Lem.2}
  If $G$ is a connected graph on $5$    vertices, then $\mathcal{H}_3(G)$ is Mengerian if and only if $G$
   is a path with double stars  or a star plus an edge    between two of its leaves.
\end{lemma}

\begin{proof}
  Sufficiency follows from Proposition  \ref{Pro.2}. Clearly, in a connected graph on $5$ vertices
  $\{v_1,v_2,v_3,v_4,v_5\}$, every vertex $v_i$ is in some connected
  subgraph on $4$ vertices and at least two vertices are not in every
  such component. If exactly two vertices are not in every such
  component, $G$ must be a path with four edges, or a path with double
  stars, or a star plus an edge between two of its leaves.  If three
  vertices, say $\{v_1,v_2,v_3\}$, are not in every path with four
  vertices, then $x_4=0, x_5=0, x_1+x_2+x_4+x_5=1= x_1+x_3= x_2+x_3=1$
  yields the solution $\frac{1}{2}(1,1,1,0,0)$, which is a vertex by
  Lemma~\ref{Lem.02}.  This implies  that the hypergraph is not ideal, and
  thus also not Mengerian. Similarly, in the case of $4$ vertices not
  in every connected component,  we find the vertex
  $\frac{1}{3}(1,1,1,1,0)$ and finally, if $G$ is $2$-connected, then  the required 
    vertex is  $\frac{1}{4}(1,1,1,1,1)$. Thus, in any of these cases,  the
  hypergraph is neither ideal nor Mengerian, and the proof is over.
\end{proof}

%%%%%%%%%%%%%%%%%%%%%%%%%%%%%%%%%%%%%%%%%%%%%%%%%%%%%%%%%%%%
%%%%%%%%%%%%%%%%%%%%%%%%%%%%%%%%%%%%%%%%%%%%%%%%%%%%%%%%%%%%

\begin{proposition} \label{Pro.3}
  If $G$ is the following graph, then $\mathcal{H}_3(G)$ is not Mengerian.
\end{proposition}

\begin{proof}
  The numbers indicated at the vertices clearly define a non-integral vertex. Thus,  $\mathcal{H}_3(G)$ is neither ideal nor Mengerian.
  
  \begin{center}
 \scalebox{1} % Change this value to rescale the drawing.
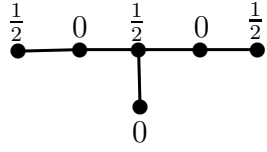
{
\begin{pspicture}(0,-0.92296875)(5.2828126,0.92296875)
\psdots[dotsize=0.2](1.0350486,0.33042014)
\psdots[dotsize=0.2](1.8550487,0.35042015)
\psdots[dotsize=0.2](2.6350486,0.35042015)
\psdots[dotsize=0.2](3.4550486,0.35042015)
\psdots[dotsize=0.2](4.215049,0.35042015)
\psdots[dotsize=0.2](2.6550486,-0.40957996)
\psline[linewidth=0.04cm](1.0550486,0.33042014)(1.8550487,0.35042015)
\psline[linewidth=0.04cm](1.8550487,0.35042015)(2.6350486,0.35042015)
\psline[linewidth=0.04cm](2.6350486,0.35042015)(3.4350486,0.35042015)
\psline[linewidth=0.04cm](3.4350486,0.35042015)(4.195049,0.35042015)
\psline[linewidth=0.04cm](2.6350486,0.37042013)(2.6550486,-0.40957996)
\usefont{T1}{ptm}{m}{n}
\rput(1.0123438,0.71453124){$\frac{1}{2}$}
\usefont{T1}{ptm}{m}{n}
\rput(4.1923437,0.7345312){$\frac{1}{2}$}
\usefont{T1}{ptm}{m}{n}
\rput(1.8523438,0.65453124){$0$}
\usefont{T1}{ptm}{m}{n}
\rput(2.6523438,-0.74546874){$0$}
\usefont{T1}{ptm}{m}{n}
\rput(2.6123438,0.71453124){$\frac{1}{2}$}
\usefont{T1}{ptm}{m}{n}
\rput(3.4723437,0.65453124){$0$}
\end{pspicture} 
}
\end{center}
\end{proof}

%%%%%%%%%%%%%%%%%%%%%%%%%%%%%%%%%%%%%%%%%%%%%%%%%%%%%%%%%%%%
%%%%%%%%%%%%%%%%%%%%%%%%%%%%%%%%%%%%%%%%%%%%%%%%%%%%%%%%%%%%

In the following result we are using a result from the  paper \cite{Candy} which has been published  hundred years ago (in 1923) and here is the first citation of this result. 

\begin{lemma} \label{Lem.3}
If $G=C_k$ is a cycle with $k\geq 5$ and $k\neq 8$, then $\mathcal{H}_3(G)$  is not ideal, and thus not Mengerian.  
\end{lemma}

\begin{proof}
  By \cite[Theorem 4.1]{CN}, these hypergraphs are not ideal, and thus
  not Mengerian either. We will present fractional vertices for
  completeness.  If $k$ is odd, then \cite[Theorem V]{Candy} implies
  that $\det(A)=4$, where $A$ denotes the incidence matrix of
  $\mathcal{H}_3(G)$.  So, $\frac{1}{4}(1, \ldots, 1)$ is a vertex. If
  $k\equiv 2$ $(\mathrm{mod}$ $4)$, then $\frac{1}{2}(1,0,1,0,1,
  \ldots, 1,0)$ is a vertex since $\frac{k}{2}$ non-negativity
  constraints and all hyperedge constraints are satisfied with
  equality and the matrix is of full rank.  Finally, if $k\equiv 0$
  $(\mathrm{mod}$ $4)$ and $k\geq 12$ using the same approach as for
  $k\equiv 0$ $(\mathrm{mod}$ $2)$ would yield a system which is not
  of full rank. But $\frac{1}{2}(1, 1, 0, 1, 0,1,1, 0, 1, 0, 1,
  \ldots, 0, 1, 0)$ is a vertex satisfying $\frac{k}{2}-1$
  non-negativity constraints and $k-4$ hyperedge constraints with
  equality. We thus have in total $\frac{3}{2}k-5 >k$ equalities, and the system is of full rank. This completes the proof.  
\end{proof}

%%%%%%%%%%%%%%%%%%%%%%%%%%%%%%%%%%%%%%%%%%%%%%%%%%%%%%%%%%%%
%%%%%%%%%%%%%%%%%%%%%%%%%%%%%%%%%%%%%%%%%%%%%%%%%%%%%%%%%%%%

A remark here what makes the case $k = 8$ special is certainly in order.

\begin{lemma} \label{Lem.4}
If $G=(V, E)$ is a connected graph with  $|V(G)|\geq 5$ and $G$ is neither $C_8$ nor a path with double stars nor  a star plus an edge  
  between two of its leaves, then $\mathcal{H}_3(G)$  is not Mengerian. 
\end{lemma}

\begin{proof}
  If $G$ is a tree which is not a path with double stars, then it has
  three vertices which are adjacent to leaves. Hence, the following  graph $H$ (see Figure 1) is a
  minor, and thus $\mathcal{H}_3(G)$ is not Mengerian by Lemma \ref{Lem.2}. If $G$ contains
  a triangle and is not itself a star plus an edge between two of its
  leaves, then we find an induced subgraph on $5$ vertices which is
  not Mengerian by Lemma~\ref{Lem.2}. 
  If $G$ contains a cycle of length $4$, again we find an induced subgraph on $5$ vertices which is not Mengerian by 
  Lemma~\ref{Lem.2}.
    
  If $G$ contains a cycle of length $k\ge 6$ and $k \ne 8$,  then  $\mathcal{H}_3(G)$ is not Mengerian by Proposition \ref{Pro.3}. 
   If the only  cycle in $G$ is the cycle  $C_8$, since $G\neq C_8$, we find a vertex
  adjacent to only one vertex of $C_8$, and so again we have the following 
   graph $H$ which is not Mengerian by Lemma \ref{Lem.2} as a minor. 
   Any other cycle is a cycle of length $\geq 6$ and the hypergraph is not Mengerian by Lemma \ref{Lem.3}.
   This finishes our argument.   
$\vspace{.2cm}$
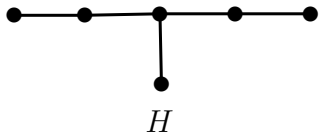
  \begin{center} \label{Figure H}
\scalebox{1} % Change this value to rescale the drawing.
{
\begin{pspicture}(0,-0.85375)(4.08,0.83375)
\psdots[dotsize=0.2](0.06,0.73375)
\psdots[dotsize=0.2](1.0,0.73375)
\psdots[dotsize=0.2](3.0,0.75375)
\psdots[dotsize=0.2](4.0,0.75375)
\psline[linewidth=0.04cm](0.06,0.75375)(0.08,0.75375)
\psline[linewidth=0.04cm](0.08,0.75375)(0.06,0.73375)
\psline[linewidth=0.04cm](0.08,0.73375)(1.04,0.73375)
\psline[linewidth=0.04cm](0.96,0.75375)(0.98,0.71375)
\psline[linewidth=0.04cm](0.98,0.69375)(1.02,0.73375)
\psline[linewidth=0.04cm](0.98,0.73375)(2.0,0.75375)
\psline[linewidth=0.04cm](2.0,0.75375)(4.0,0.75375)
\usefont{T1}{ptm}{m}{n}
\rput(1.9914062,-0.67625){$\text{Figure 1. The graph}~ H$}
\psdots[dotsize=0.2](2.0,0.75375)
\psline[linewidth=0.04cm](1.98,0.77375)(2.04,0.73375)
\psline[linewidth=0.04cm](2.0,0.73375)(2.02,-0.22625)
\psdots[dotsize=0.2](2.02,-0.18625)
\end{pspicture} 
}
  \end{center}
 \end{proof}

%%%%%%%%%%%%%%%%%%%%%%%%%%%%%%%%%%%%%%%%%%%%%%%%%%%%%%%%%%%%
%%%%%%%%%%%%%%%%%%%%%%%%%%%%%%%%%%%%%%%%%%%%%%%%%%%%%%%%%%%%

The main result of this paper follows combining Proposition  \ref{Pro.2} with Lemma \ref{Lem.4}.

\begin{theorem} \label{Main-Result}
If $G=(V, E)$ is a connected graph, then $\mathcal{H}_3(G)$ is Mengerian if and only if $|V|=4$, or $G=C_8$, or $G$ is a path with double stars, or 
a star with an additional edge  between two of its leaves.
\end{theorem}

\begin{corollary}
If $G=(V, E)$ is a connected graph with  $|V|=4$, or $G=C_8$, or $G$ is a path with double stars, or 
a star with an additional edge  between two of its leaves, then $\mathcal{H}_3(G)$ satisfies the Conforti-Cornu\'e{j}ols  conjecture. 
\end{corollary}

%%%%%%%%%%%%%%%%%%%%%%%%%%%%%%%%%%%%%%%%%%%%%%%%%%%%%%%%%%%%
%%%%%%%%%%%%%%%%%%%%%%%%%%%%%%%%%%%%%%%%%%%%%%%%%%%%%%%%%%%%

\section{Conclusion and Outlook} 
In this paper we considered the
4-uniform hypergraphs $\mathcal{H}_3(G)$, where the hyperedges are
defined by the vertex sets $h \in \binom{V}{4}$ of a graph $G=(V,E)$,
that contain a path on $3$ edges.  We characterized the graphs where
the corresponding hypergraphs are Mengerian. Except for
$\mathcal{H}_3(C_8)$, the corresponding matrices where all even
totally unimodular. Moreover, we found a dichotomie that in this
class every hypergraph is either Mengerian or the corresponding matrix
is non-ideal. In view of Theorem 4.1(i) in~\cite{CN}, we wonder whether

\begin{question} \label{Open Question}
    Is it true that for $t\ge 4$ every matrix of the $(t+1)$-uniform
  hypergraph $\mathcal{H}_t(G)$ either is totally unimodular or
  non-ideal?
\end{question}

%%%%%%%%%%%%%%%%%%%%%%%%%%%%%%%%%%%%%%%%%%%%%%%%%%%%%%%%%%%%
%%%%%%%%%%%%%%%%%%%%%%%%%%%%%%%%%%%%%%%%%%%%%%%%%%%%%%%%%%%%

\section*{ORCID}

  Winfried  Hochstättler:  \textsc{https://orcid.org/0000-0001-7344-7143} \par
Mehrdad ~Nasernejad: \textsc{https://orcid.org/0000-0003-1073-1934}\\

%%%%%%%%%%%%%%%%%%%%%%%%%%%%%%%%%%%%%%%%%%%%%%%%%%%%%%%%%%%%
%%%%%%%%%%%%%%%%%%%%%%%%%%%%%%%%%%%%%%%%%%%%%%%%%%%%%%%%%%%%
 
\section*{The conflict of interest and data availability statement}

We hereby declare  that this  manuscript has no associated  data and also there is no conflict of interest in this manuscript.

%%%%%%%%%%%%%%%%%%%%%%%%%%%%%%%%%%%%%%%%%%%%%%%%%%%%%%%%%%%%
%%%%%%%%%%%%%%%%%%%%%%%%%%%%%%%%%%%%%%%%%%%%%%%%%%%%%%%%%%%%

\section*{Acknowledgment.} 
First, the authors  are  deeply grateful to the anonymous referee for careful reading of the manuscript, and for  his/her
  valuable suggestions which led to several  improvements in the quality of this paper and pointed out a gap in Lemma~\ref{Lem.01}. 
In addition, this paper was prepared when the second author  visited the Department of Mathematics of 
the University of Hagen  in 2023. In particular, he  would like to thank the University of Hagen  for its hospitality.
Furthermore, this work has been supported by   DAAD (Germany) (Funding  No.  57681226).

%%%%%%%%%%%%%%%%%%%%%%%%%%%%%%%%%%%%%%%%%%%%%%%%%%%%%%%%%%%%
%%%%%%%%%%%%%%%%%%%%%%%%%%%%%%%%%%%%%%%%%%%%%%%%%%%%%%%%%%%%
%%%%%%%%%%%%%%%%%%%%%%%%%%%%%%%%%%%%%%%%%%%%%%%%%%%%%%%%%%%%
%%%%%%%%%%%%%%%%%%%%%%%%%%%%%%%%%%%%%%%%%%%%%%%%%%%%%%%%%%%%

%%%%%%%%%%%%%%%%%%%%%%%%%%%%%%%%%%%%%%%%%%%%%%%%%%%%%%
%%%%%%%%%%%%%%%%%%%%%%%%%%%%%%%%%%%%%%%%%%%%%%%%%%%%%%

% ------------------------------------------------------------------------
\end{document}